%% file: main.tex
\documentclass[twocolumn,final]{autart}
\usepackage{graphicx}
\usepackage{xcolor}
\usepackage{amsmath,amsfonts}
\usepackage{mathtools,cuted}
\usepackage{subfiles}
\usepackage{subcaption}
\usepackage{enumitem}
\usepackage{todonotes}
\usepackage[numbers,sort&compress]{natbib}
\usepackage{tikz}
\usetikzlibrary{decorations.markings}
\usetikzlibrary{shapes,arrows,positioning}
\usetikzlibrary{arrows.meta}
\usepackage[siunitx, americanvoltages, americancurrents,smartlabels]{circuitikz}
\tikzstyle{every node}=[font=\small]
\usepackage{mdframed}
\usetikzlibrary{matrix}
\usepackage{soul}
\usepackage[normalem]{ulem}
\usepackage{cancel}

\newcommand{\ve}{\varepsilon}
\newcommand{\R}{\mathbb{R}}

\newcommand{\bx}{\bar{x}}
\newcommand{\bz}{\bar{z}}

\newcommand{\br}{\bar{r}}
\newcommand{\be}{\bar{\ve}}
\newcommand{\oX}{\overline{X}}
\newcommand{\tX}{\widetilde{X}}

\newcommand{\parcs}[2]{\frac{\partial #1}{\partial #2}}

\DeclareMathOperator{\diag}{diag}

\newtheorem{remark}{Remark}
\newtheorem{definition}{Definition}
\newtheorem{theorem}{Theorem}

\newcommand{\change}[1]{\color{black}{#1}}

\begin{document}

\begin{frontmatter}
	\title{Stabilization of slow-fast control systems: the non-hyperbolic case}
	\author[rug]{Hildeberto Jard\'on-Kojakhmetov}\ead{h.jardon.kojakhmetov@rug.nl},    
	\author[rug]{Jacquelien M.A. Scherpen}\ead{j.m.a.scherpen@rug.nl},

	\address[rug]{Jan C. Willems Center for Systems and Control, Faculty of Science and Engineering (FSE), Engineering and Technology Institute (ENTEG), University of Groningen, The Netherlands}


	\begin{abstract}
	In this paper we study the stabilization problem of a general class of slow-fast systems with one fast and arbitrarily many slow states. Moreover, the class of systems under study is \emph{slowly actuated}, meaning that only the slow states are subject to the action of a controller. Furthermore, we are particularly interested in the case where normal hyperbolicity is lost. We show that by using the Geometric Desingularization method, it is possible to design controllers to locally stabilize non-hyperbolic points of any finite degeneracy. The main novelty of this paper is that, unlike previous research on the topic, we make use of more than one chart of the blow up space to enhance the region of attraction of the operating point. A couple of numerical examples highlight our contribution.
\end{abstract}

\end{frontmatter}

\section{Introduction}
\subfile{subfiles/intro.tex}

\section{Preliminaries}\label{sec:preliminaries}
\subfile{subfiles/pre.tex}

\section{Geometric Desingularization}\label{sec:gd}
\subfile{subfiles/geom.tex}

\section{Main results}\label{sec:results}
\subfile{subfiles/main-res.tex}

\section{Examples}\label{sec:sim}

\subfile{subfiles/ex.tex}

\section{Conclusions}\label{sec:conclusions}
In this paper we have shown a novel method to design controllers for slow-fast control systems that render a non-hyperbolic point asymptotically stable. To this end we have used a technique called Geometric Desingularization together with simple control techniques.  In particular, we have provided a new control methodology for slow-fast control systems that classical techniques \cite{Kokotovic:1986:SPM:576779} do not cover. As future research we propose the extension of the presented methodology to trajectory and path following control problems. Special difficulties arise when the trajectory or path to be followed passes through non-hyperbolic points. Another complicated problem is to study slow-fast control systems near non-hyperbolic points with more than one fast direction.





\bibliographystyle{plain}
\bibliography{Mendeley}
\end{document}

%% file: subfiles/intro.tex
Multiple timescales are ubiquitous in mathematical modeling and applications. Examples of real life phenomena with several timescales can be found in nonlinear circuits \cite{vanderpol_heart,smale1972mathematical,ihrig1975regularization,Reissig1,RIAZA2011}, neuron models \cite{Shilnikov2012}, biochemical systems \cite{Rotstein2013,kosiuk2016geometric}, etc. There exists a great body of literature dealing with systems for which the timescale separation is global \cite{Fenichel1979,Kokotovic:1986:SPM:576779}. However, in many complex models, the aforementioned global timescale separation does not hold. Usually, from a dynamical systems point of view, a non-global timescale separation is associated with the presence of certain singularities. The correct and thorough analysis of the behavior of a multi-timescale system near such singularities is crucial for the progress of our understanding of phenomena with several timescales. 

In the context of control systems, considerable effort has been given to study systems with two timescales and global separation of such timescales \cite{ Kokotovic:1986:SPM:576779,KokotovicApps,Marino1988,Valery,Saksena1984}. However, with the increased interest in shaping the behavior of complex multi-timescale systems, we require control techniques that can tackle problems where the timescale separation is not global. Preliminary steps in this regard have been developed for regulation purposes in \cite{JardonMTNS2016,JardonACC2017} in the planar case, and \cite{jardon2017stabilization} for the case of singularities of quadratic degeneracy. 

In this {\change{article}} we extend the results of \cite{jardon2017stabilization} to a  broader class of slow-fast control systems, and we propose a geometric way to enlarge the region of stability of an equilibrium point. The idea is to show that Geometric Desingularization \cite{DumRou2,Kuehn2015}  can be used in combination with control strategies to stabilize degenerate points of slow-fast systems, c.f. \cite{JardonMTNS2016,JardonACC2017} for the planar case. In few words, the Geometric Desingularization technique is a suitable change of coordinates, well-defined around singularities, that allows a simpler analysis of the involved behavior of slow-fast systems without global timescale separation. In turn, by employing such a technique, the design of controllers for slow-fast systems in the aforementioned scenario becomes simpler. 

{\change{Briefly speaking, we study slow-fast control systems with one fast and an arbitrary amount of slow states. Moreover, the system us actuated \emph{only} on the slow variables. Our contributions can be summarized as follows: first, we show that by using the Geometric Desingularization technique we can design, in a rather simple way, controllers that locally stabilize a non-hyperbolic point of a slowly actuated slow-fast systems as described above (Theorem \ref{thm:main}). Next, we provide a constructive way to design a controller that accomplishes the aforementioned task (Theorem \ref{thm:main2}). Furthermore, the main novelty of this article is presented in Theorem \ref{thm:reg} where we show that with a more thorough analysis of the desingularized closed-loop system, we can enlarge the region of attraction of an equilibrium point. Finally, as a benchmark, we compare our proposed controller with a high-gain one, and show that our proposed controller requires a much smaller gain than the high-gain controller in order to stabilize a non-hyperbolic point of a slow-fast control system .}}

\change{The rest of this document is organized as follows: in Section \ref{sec:preliminaries} we provide the necessary preliminaries and the formal setting of the problem to be studied. Next, in Section \ref{sec:gd} we briefly recall the Geometric Desingularization technique, which is essential to prove our results. Afterwards, in Section \ref{sec:results}, we present our contributions, followed by two examples in Section \ref{sec:sim} highlighting our results. Finally, Section \ref{sec:conclusions} presents some concluding remarks and a couple of open problems for future research. }

%% file: subfiles/pre.tex
	From now on, we shall confine ourselves to the study of two-timescale systems, also known as slow-fast system (SFS). A SFS is an Ordinary Differential Equation of the form
\begin{equation}\label{sfs0}
\begin{split}
\dot x &= f(x,z,\ve)\\
\ve\dot z &= g(x,z,\ve),
\end{split}
\end{equation}
where $x\in\R^{n_s}$ are the slow states, $z\in\R^{n_f}$ the fast states, $0<\ve\ll1$ is a small parameter responsible for the timescale separation between $x$ and $z$, and $f$ and $g$ are smooth functions. For $\ve>0$ we can define a new time parameter $\tau=t/\ve$ and obtain an equivalent system to \eqref{sfs0} of the form
\begin{equation}\label{sfs00}
\begin{split}
x' &= \ve f(x,z,\ve)\\
z' &= g(x,z,\ve).
\end{split}
\end{equation}%
Usually, in the study of SFSs, we define two reduced subsystems by taking the limit $\ve\to 0$ of \eqref{sfs0} and \eqref{sfs00}. By doing so we obtain the DAE (Differential Algebraic Equation) and the Layer Equation, which read as
\begin{minipage}{0.25\textwidth}
\begin{equation*}\label{cde}
\text{DAE}:\begin{cases}
\dot x &= f(x,z,0)\\
0 &= g(x,z,0)
\end{cases}
\end{equation*}
\end{minipage}%
\begin{minipage}{0.25\textwidth}
\begin{equation*}\label{lay}
\text{Layer}:\begin{cases}
x' &= 0\\
z' &= g(x,z,0)
\end{cases}
\end{equation*}
\end{minipage}
The main idea of (Geometric) Singular Perturbation Theory is to draw conclusions (for example: qualitative and quantitative description, stability, etc.) from the reduced systems {\change{(i.e. the DAE and the Layer Equation)}} and then extend them to similar results of the corresponding SFS. In the analysis of SFSs, one of the most important geometric objects to consider is the \emph{critical manifold}.
\begin{definition} The critical manifold associated to the SFS \eqref{sfs0} is defined as
\begin{equation}
S=\left\{ (x,z)\in\R^{n_s+n_f}\, | \, g(x,z,0)=0 \right\}.
\end{equation}
\end{definition}
Note that $S$ is the phase-space of the DAE and the set of equilibrium points of the Layer Equation. We say that $S$ is \emph{Normally Hyperbolic} if every point $s\in S$ is a hyperbolic equilibrium point of the reduced dynamics $z'=g(x,z,0)$. The theory regarding SFSs with Normally Hyperbolic critical manifold is nowadays well understood, see e.g. \cite{Fenichel1979} for a general treatment, and \cite{Kokotovic:1986:SPM:576779} for applications in the context of control systems. However, we still find many open problems in the situation where $S$ has non-hyperbolic points. 

\subsection{Setting}

The class of slow-fast control systems (SFCSs) that we study are defined as 
\begin{equation}\label{sfs1}
\begin{split}
\dot x &= f(x,z,u,\ve)\\
\ve\dot z &= g(x,z,\ve),
\end{split}
\end{equation}
where {\change{$x\in\R^{n_s}$, $n_s\geq 1$ is an integer}}, $z\in\R$, $f(x,z,0,0)$ is smooth, {\change{$u\in\R^{n_s}$}} denotes the control input, and $g(x,z,\ve)$ is a smooth function satisfying
\begin{equation}\label{g}
g(x,z,0)= -\left( z^k + \sum_{i=1}^{k-1} x_iz^{i-1} \right),
\end{equation}%
{\change{where $k\in\mathbb N$, with $k\geq 2$.}}
\begin{remark}\leavevmode
\begin{itemize}[leftmargin=*]
\item The control input $u$ only acts on the slow dynamics, that is, the class of systems \eqref{sfs1}-\eqref{g} is under-actuated.
\item The origin is the most degenerate non-hyperbolic point of the critical manifold $S=\left\{ g(x,z,0)=0 \right\}$. {\change{ This is because the origin is the unique point where $g(0)=\cdots=\frac{\partial^{k-1} g}{\partial z^{k-1}}(0)=0$, and $\frac{\partial^{k} g}{\partial z^{k}}(0)\neq0$ }} Therefore, we are interested in stabilizing the origin of \eqref{sfs1}.
\item The fact that the class of systems \eqref{sfs1}-\eqref{g} is large is explained by the classification of singularities of smooth functions. Let us give a brief recollection of the relevant arguments, for more details see \cite{arnold1974critical,Arnold_singularities,Brocker}. Let $V_0(z)$, where $z\in\R$, be a smooth function with $V_0(0)=0$. Let $k\geq 0$ be the minimum integer such that $\frac{\partial^0 V_0}{\partial^0 z}(0)= \cdots=\frac{\partial^k V_0}{\partial^k z}(0)= 0$ and $\frac{\partial^{k+1} V_0}{\partial^{k+1} z}(0)\neq 0$. Then $V_0(z)$ is locally equivalent to $\pm z^{k+1}$ (see e.g. \cite[Theorem 1]{arnold1974critical}). 
Next, let $V_x(z)=V(x,z)$ be a generic family of functions with $V_0(z)$ as above. Then, the universal unfolding (containing the minimum number of parameters such that the singularity is generic) of $V_0(z)$ is $V_x(z)=\pm z^{k+1} + x_{k-1}z^{k-1}+\cdots+x_1z$ (see e.g. \cite[Example 14.9]{Brocker}). Thus, the critical manifold $S$ associated to \eqref{sfs1} is equivalently defined as the set of critical points of $V_x(z)$. {\change{The fact that we need at least $k-1$ parameters to unfold $V_0$ implies that it suffices to fix $n_s=k-1$. If $n_s<k-1$, then the class of SFCSs \eqref{sfs1}-\eqref{g} is not generic, on the other hand, if $n_s>k-1$ similar techniques as used here can be employed, compare with \cite{jardon2017stabilization}}}. The choice of the negative sign in \eqref{g} is just for convenience. {\change{It simply means that away from $S$, the trajectories of the layer equation travel ``downwards''. Choosing a positive sign in  \eqref{g} reverses the aforementioned direction, but a similar analysis as the one performed here can be used in such a case}}. 
\end{itemize}
\end{remark}

%% file: subfiles/geom.tex
\documentclass[../main.tex]{subfiles}
The Geometric Desingularization method is used to overcome the difficulties that the presence of non-hyperbolic points pose. Briefly speaking, this method provides a new system,  equivalent to \eqref{sfs1}, but with simpler singularities. In this way, the design of controllers for \eqref{sfs1} becomes more accessible. The Geometric Desingularization method is motivated by the regularization of singularities in algebraic varieties \cite{Hironaka1,Hironaka2}. In the context of slow-fast systems it was first introduced in \cite{DumRou2}, and has been further developed afterwards, see e.g. \cite{Jardon-Kojakhmetov2016AnalysisSingularity,Jardon-Kojakhmetov2015,Krupa1,Krupa20102841,Kuehn2015} and \cite{jardon2017stabilization,JardonMTNS2016,JardonACC2017,jardon2017modelred} for some applications in control systems.

To start the description of the method, let us rewrite a SFCS as
\begin{equation}\label{eq:sfcs}
X=\ve f(x,z,u,\ve)\parcs{}{x}+g(x,z,\ve)\parcs{}{z}+0\parcs{}{\ve},
\end{equation}
which is a smooth vector field on $\R^{k+1}$. The Geometric Desingularization method consists of the following steps (see more details in \cite{jardon2017stabilization,Kuehn2015}).
\begin{enumerate}[leftmargin=*]
\item Define the quasi-homogeneous blow up\footnote{In the context of multi-timescale systems, the term blow up is understood as a ``zoom-in'' and not as an explosion.} map $\Phi:\mathbb S^k\times I\to\R^{k+1}$, where $I\subseteq [0,\infty)$ is an interval  (possibly infinite), and $\mathbb S^k$ denotes the $k$-sphere, by
\begin{equation}\label{eq:bu}
\begin{split}
\Phi(\bx,\bz,\be,r) &=(r^\alpha\bx,r^\beta\bz,r^\gamma\be)=(x,z,\ve),
\end{split}
\end{equation}
where $\bx=(\bx_1,\ldots,\bx_{k-1})$, $(\bx,\bz,\be)\in\mathbb S^k$, that is $\sum_{i=1}^{k-1}\bx_i^2+\bz^2+\be^2=1$, 
$r^\alpha\bx=(r^{\alpha_1}\bx_1,\ldots,r^{\alpha_{k-1}}\bx_{k-1})$, $\alpha=(\alpha_1,\ldots,\alpha_{k-1})\in\mathbb N^{k-1}$, $\beta\in\mathbb N$, $\gamma\in\mathbb N$, and $r\in I$. If the weights are all equal to $1$, we simply refer to $\Phi$ as blow up.

\item Define the desingularized vector field $\tX$ as follows: i) Note that the blow up map $\Phi$ induces a vector field $\oX$ on $\mathbb S^k\times I$ defined as $\oX=D\Phi^{-1}\circ X\circ\Phi$, where $D\Phi$ denotes the directional derivative of $\Phi$. Since $\Phi$ is a diffeomorphism for $\left\{ r>0 \right\}$, the vector fields $X$ and $\oX$ are conjugate for all $r>0$. Moreover, the definition of $\oX$ extends continuously to $\left\{ r=0\right\}$ \cite{Kuehn2015}, that is $\oX$ is well defined on $\mathbb S^k\times I$. Note that $\mathbb S^k\times \left\{ 0 \right\}$ is mapped to the origin $0\in\R^{k+1}$ via the blow up map, therefore, since $X(0)=0$, $\oX$ vanishes along $\mathbb S^k\times \left\{ 0 \right\}$, so ii) Define the \emph{desingularized} vector field $\tX$ as $\tX=\frac{1}{r^p}\oX$, where $p\in\mathbb N$ is suitably chosen so that $\tX$ does not vanish along $\mathbb S^k\times \left\{ 0 \right\}$.
\end{enumerate}

Note that the vector fields $\oX$ and $\tX$ are smoothly equivalent on $\mathbb S^k\times\left\{ r>0 \right\}$, their only difference is the time-parametrization. Thus, it is qualitatively the same to study $\tX$ than $\oX$. This also implies that the qualitative properties of $\tX$ can be related to similar ones of $X$, after all, they are also smoothly equivalent for all $r>0$. The idea of Geometric Desingularization is to appropriately choose the weights of the blow up map \eqref{eq:bu} so that the vector field $\tX$ has simpler singularities (e.g. hyperbolic, or semi-hyperbolic) with respect to those of the SFCS  \eqref{eq:sfcs}, and that the singular dependence of $X$ on $\ve$ is overcome in the blow up space. In turn, the design of the controller for $\tX$ becomes simpler.

\begin{remark}
Usually, if the vector field $X$ is quasihomogeneous \cite{jardon2017stabilization,Kuehn2015} then the weights of the blow up map $\Phi$ correspond to the quasihomogeneity type of $X$, but in general, finding the appropriate values of the weights is a non-trivial task. However, as every singularity of an algebraic variety over a field of characteristic $0$ (such as $\R$) can be regularized after a finite number of blow ups \cite{Hironaka1,Hironaka2}, we {\change{conjecture}} that, for any vector field on $\R^n$, this methodology also holds.
\end{remark}

While performing the computations, it is more convenient to introduce charts rather than working on spherical coordinates. A chart is just a parametrization of a hemisphere of the blown up space. In practice, a chart is obtained by simply setting one of the coordinates $(\bx_1,\ldots,\bx_{k-1},\bz,\be)\in\mathbb S^k$ to $\pm1$ in the definition of $\Phi$. In this way we define, for example, the chart $K_{\be}=(r^\alpha\bx,r^\beta\bz,r^\gamma)$. Note that each chart parametrizes just a part of $\mathbb S^k\times I$; however, all possible charts define an open cover of $\mathbb S^k\times I$. To avoid confusion, whenever we work on more than one chart, we shall define local coordinates and distinguish them in each chart. The charts are related to each other via \emph{transition maps} \cite{jardon2017stabilization,Kuehn2015}.

\begin{remark}
The chart $K_{\be}$ is the most important one in our analysis, and it is called \emph{the family chart}. All other charts, which we denote by $K_{\pm\bx_i}=\left\{ \bx_i=\pm1\right\}$ and $K_{\pm\bz}=\left\{ \bz=\pm1\right\}$ are called \emph{directional charts}. In the chart $K_{\be}$ the singular dependence of the vector field on $\ve$ is overcome. It is also the chart where most of the local dynamical properties of $X$ near the origin can be seen.
\end{remark}

%% file: subfiles/main-res.tex
Here we present our main results: we show that a controller designed for the desingularized system $\tX$ gives, after change of coordinates, a controller for $X$. 

\begin{theorem}\label{thm:main} Consider the SFCS $X$ defined by \eqref{eq:sfcs} and suppose that $X$ can be desingularized as described in Section \ref{sec:gd}. Denote the desingularization of $X$ by $\tX$. If $\widetilde u$ is a controller that renders $\mathbb S^k\times\left\{ 0 \right\}$ asymptotically stable\footnote{ We say that the set $B_0=\mathbb S^k\times\left\{ 0 \right\}$ is asymptotically stable if there is a neighborhood $\mathcal N$ of $B_0$ such that every trajectory $\bar \gamma(t)$ of $\tX$ with initial condition in $\mathcal N$ satisfies $d_H(\gamma(t),B_0)\to 0$ as $t\to\infty$, where $d_H$ denotes Hausdorff distance. }, under the flow of $\tX$, then $u(x,z,\ve)=\Phi\circ\widetilde u\circ\Phi^{-1}(x,z,\ve)$ renders the origin $0\in\R^{k+1}$ asymptotically stable under the flow of $X$.		
\end{theorem}
\begin{pf} Let us rewrite the closed-loop blown up system as $\widetilde Y$, that is $\widetilde Y(\bx,\bz,\be,r)=\tX(\bx,\bz,\be,r,\widetilde u)$, where $\widetilde u=\widetilde u(\bx,\bz,\be,r)$ is a feedback controller that renders $B_0=\mathbb S^k\times\left\{ 0 \right\}$  asymptotically stable. Similarly we denote by $Y$ the closed-loop system $X$ with $u$ induced by $\Phi$, that is $u(x,z,\ve)=\Phi\circ\widetilde u\circ\Phi^{-1}(x,z,\ve)$. Therefore $\widetilde Y$ is the desingularization of $Y$. Recall that $\Phi$ maps $\mathbb S^k\times\left\{ 0 \right\}$ to $0\in\R^k$ and that $Y$ and $\widetilde Y$ are equivalent. This means that trajectories of $\widetilde Y$ are mapped, via $\Phi$, to trajectories of $Y$ preserving direction. Therefore, trajectories of $\widetilde Y$ approaching $\mathbb S^k\times\left\{ 0 \right\}$ are mapped to trajectories of $Y$ that approach $0\in\R^k$. The asymptotic convergence is preserved since $\Phi$ is a diffemorphism. Note that by taking $r\in I$ with $I$ arbitrarily large, we can find an equivalence between the trajectories of $Y$ and of  $\widetilde Y$ within an arbitrarily large compact subset of $\R^k$. \hfill\ensuremath\qed
\end{pf}

The main idea of Theorem \ref{thm:main} is that we can design controllers for a SFCS $X$ in the blow up space. Due to the fact that the blow up simplifies the singularities of a vector field, the design of controllers in the blow up space is expected to be simpler than in the original scenario, i.e., without blow up. Note, however, that Theorem \ref{thm:main} requires an analysis within the whole blow up space. This may be computationally tedious, so we also present a more relaxed result and with a particular choice of controller, which for certain applications may be sufficient. The main idea of the following result is to design a feedback controller in the blow up space that linearizes the slow dynamics.

\begin{theorem}\label{thm:main2} Consider a SFCS given by
\begin{equation}\label{eq:th2}
\begin{split}
\dot x &= f(x,z,\ve) +u\\
\ve\dot z &= -\left( z^k+\sum_{i=1}^{k-1}x_i z^{i-1}\right),
\end{split}
\end{equation}
where $x=(x_1,\ldots,x_{k-1})\in\R^{k-1}$, $z\in\R$ and $u\in\R^{k-1}$. 
Then, for $\ve>0$ sufficiently small, the controller
\begin{equation}\label{eq:u-thm2}
u=-C+b\ve^{\frac{-1}{2k-1}}z\hat e_1-\ve^{\frac{-k}{2k-1}}Ax,
\end{equation}
where $C=\diag\left\{f_i(0,0,0)\right\}$, $b>0$, $A>0$ diagonal, and $\hat e_1=[1,\, 0, \, \cdots \,, 0]^\top\in\R^{k-1}$, renders the origin $(x,z)=(0,0)\in\R^k$ locally asymptotically stable.
\end{theorem}
\begin{pf}
	The proof consists of designing the controller in the chart $K_{\bar\ve}$. {\change{Therefore, the very first step on this proof is to find appropriate weights for the blow up map $\Phi$. In what follows, we show one possible way to choose such weights.}} Let us first rewrite \eqref{eq:th2} as the vector field
\begin{equation}
X=\ve\left( f(x,z,\ve)+u\right)\parcs{}{x}+g(x,z)\parcs{}{z}+0\parcs{}{\ve}.
\end{equation}

Note that $g(x,z)=-\left( z^k+\sum_{i=1}^{k-1}x_i z^{i-1}\right)$ is a quasihomogeneous polynomial of type $(k,k-1,\ldots,1)$. Thus, let us propose the coordinate transformation
\begin{equation}
(x_1,\ldots,x_{k-1},z,\ve)=(\br^k\bx_1,\ldots,\br^2\bx_{k-1},\br\bz,\br^\gamma),
\end{equation}
where $\gamma\change{\in\mathbb N}$ shall be {\change{appropriately chosen}} below. Next, we obtain the blown up vector field. First, from $\ve'=0$ and $\ve=\br^\gamma$, it follows that $\br'=0$. Similarly we obtain
\begin{equation}\label{eq:pp1}
\begin{split}
\bx_i' &= \br^{\gamma-k+i-1}(\bar f_i +\bar u _1), \qquad i=1,2,\ldots,k-1\\
\bz' &= -\br^{k-1}\left( \bz^k+\sum_{i=1}^{k-1}\bx_i \bz^{i-1}\right),
\end{split}
\end{equation}
where $\bar f_i=\bar f_i(\br,\bx,\bz)=f(\br^k\bx_1,\ldots,\br^2\bx_{k-1},\br\bz,\br^{2k-1})$, and similar notation is used for $\bar u_i$, $i=1,\ldots,k-1$. Then, to desingularize \eqref{eq:pp1}, we need to divide the right hand side by $\br^{k-1}$. {\change{To obtain a well-defined desingularized vector field when $\br=0$, it is convenient to choose $\gamma=2k-1$. A choice such that $\gamma<2k-1$ does not provide a well defined vector field for the restriction $r=0$. On the other hand, the choice $\gamma>2k-1$} would imply $\bx_i'\in O(\br)$ for all $i=1,\ldots,k-1$.}   With this we obtain the desingularized system
\begin{equation}\label{eq:pp2}
\begin{split}
\bx_i' & = \br^{i-1}(\bar f_i + \bar u_i), \qquad i=1,2,\ldots,k-1\\
\bz' &= -\left( \bz^k+\sum_{i=1}^{k-1}\bx_i \bz^{i-1}\right).
\end{split}
\end{equation}

Next, since we can pick $\bar u_i$ arbitrarily, let us propose that each $\bar u_i$ is of the form  $\bar u_1=-c_1-a_1\bx_1+b\bz$, and $\bar u_i=-c_i-\br^{1-i}a_i\bx_i$ for $i=2,\ldots,k-1$ and with $c_i=f_i(0,0,0)$, $a_i>0$, and $b>0$. Then \eqref{eq:pp2} is rewritten as
\begin{equation}\label{eq:pp3}
\begin{split}
\bx_1' &=\bar f_1 -c_1-a_1\bx_1+b\bz\\
\bx_2' &= \br( \bar f_2-c_2)-a_2\bx_1\\
&\vdots \\
\bx_{k-1}' &= \br^{k-2}( \bar f_{k-1}-c_{k-1})-a_{k-1}\bx_{k-1}\\
\bz' &= -\left( \bz^k+\sum_{i=1}^{k-1}\bx_i \bz^{i-1}\right).
\end{split}
\end{equation}   
The role of $c_i$ is to eliminate the constant values of $f_i(0,0,0)$, while $a_i$ and $b$ are chosen so that the origin is locally asymptotically stable. Note that \eqref{eq:pp3} depends regularly on $\br$. Therefore, it is convenient to study the stability of \eqref{eq:pp3} for $\br=0$. Then, the same stability properties hold for $\br>0$ sufficiently small. We remark at this point that from the relation $\ve=\br^{2k-1}$, the arguments regarding the stability of \eqref{eq:pp3} for $\br>0$ sufficiently small are equivalent to similar ones for \eqref{eq:th2} with $\ve=\br^{1/(2k-1)} $ sufficiently small. Note that $\bar f_i(0,\bx,\bz)=f_i(0,0,0)=c_i$. Therefore, system \eqref{eq:pp3} restricted to $\br=0$ reads as
\begin{equation}\label{eq:pp4}
\begin{split}
\bx_1' &= -a_1\bx_1+b\bz\\
\bx_2' &= -a_2\bx_1\\
&\vdots \\
\bx_{k-1}' &= -a_{k-1}\bx_{k-1}\\
\bz' &= -\left( \bz^k+\sum_{i=1}^{k-1}\bx_i \bz^{i-1}\right).
\end{split}
\end{equation}
The Jacobian of \eqref{eq:pp4} evaluated at the origin is of the form
\begin{equation}
J=\begin{bmatrix}
-A & b\hat e_1\\ 
-\hat e_1^\top & 0
\end{bmatrix},
\end{equation}
where $e_1=(1,0,\ldots,0)^\top\in\R^{k-1}$ and $A$ is a $(k-1)\times(k-1)$ diagonal matrix of the form $A=\diag\left\{ a_i\right\}$. It is straightforward to show that the eigenvalues of $J$ are 
\begin{equation}
\left\{ \frac{-a_1\pm\sqrt{a_1^2-4b}}{2},-a_2,\ldots,-a_{k-1} \right\}.
\end{equation}
Thus, we conclude that the origin of \eqref{eq:pp4} is locally asymptotically stable. To obtain the controller $u$, we just perform the blow down (inverse of the blow up) resulting in $u=(u_1,\ldots,u_{k-1})$ with
\begin{equation}
\begin{split}
u_1 &= -c_1-a_1\ve^{\frac{-k}{2k-1}}x_1-b\ve^{\frac{-1}{2k-1}}z\\
u_j &= -c_2-a_j\ve^{\frac{-k}{2k-1}}x_j, \qquad j=2,\ldots,k-1,
\end{split}
\end{equation}
which in compact form can be presented as in \eqref{eq:u-thm2}. \hfill\ensuremath\qed
\end{pf}

The proof of Theorem \ref{thm:main2} consists of designing a controller in the family chart $K_{\bar\ve}$. A limitation of Theorem \ref{thm:main2} is that, upon choosing $A$ and $b$, the region of attraction of the origin shrinks as $\ve\to 0$. This can be resolved by studying the system in the directional  charts to upgrade and improve the controller, as shown in the following section.

\subsection{On the region of attraction}

We now discuss a way to extend the region of attraction of the origin via further control actions. We remark that we do this without any Lyapunov analysis but rather by studying the dynamics of the blown up system in an appropriate chart. Increasing the region of attraction is particularly important for SFCS \eqref{eq:sfcs} with $k$ even. This is because, {\change{near the origin, the critical manifold is "$U$-shaped" when $k$ is even while it is "$S$-shaped" when $k$ is odd. This means that, for $k$ even, trajectories that do not converge to the origin quickly diverge. See Figure \ref{foldandcusp} for a schematic impression of the previous description }}. 
{\change{
\begin{figure}[htbp]\centering
\begin{minipage}{0.25\textwidth}
\begin{tikzpicture}
\pgftext{\includegraphics[scale=0.9]{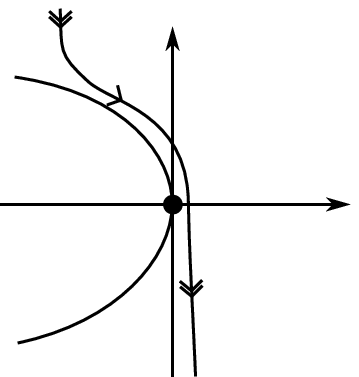}}
\node at (1.75,-.1) {$x$};
\node at (0,1.75) {$z$};
\end{tikzpicture}
\end{minipage}%
\begin{minipage}{0.25\textwidth}
\begin{tikzpicture}
\pgftext{\includegraphics[scale=0.9]{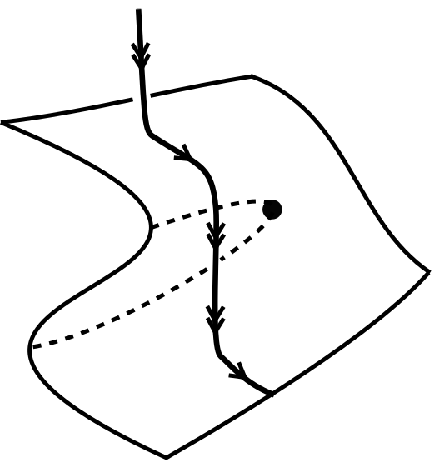}}
\end{tikzpicture}
\end{minipage}%
\caption{\change{ {\textbf{Left:}} Example of a critical manifold for the case $k$ even, given by $S=\left\{ (x,z)\in\R^2\,|\, z^2+x=0 \right\}$. The origin is a fold point. {\textbf{Right:}} Example of a critical manifold for the case $k$ odd, given by $S=\left\{ (x_1,x_2,z)\in\R^3\,|\, z^3+x_2z+x_1=0 \right\}$. The origin is a cusp point, while the dashed line represents fold points. Note that the trajectory jumps when it passes through a fold point and not through a cusp point. }}
\label{foldandcusp}
\end{figure}
}}



\begin{theorem}\label{thm:reg}
Consider the SFCS \eqref{eq:th2}. Suppose that the controller is given by
\begin{equation}\label{eq:cuu}
u=-C+b\ve^{\frac{-1}{2k-1}}z\hat e_1-\ve^{\frac{-k}{2k-1}}Ax+w,
\end{equation}
where $C$, $b$, and $A$ are as in Theorem \ref{thm:main2}, and $w=(w_1,\ldots,w_{k-1})\in\R^{k-1}$ is given by
\begin{equation}\label{eq:cw}
w_i=K_i\left(x_iz+(-z)^{k-i+2}\chi_i^*\right),
\end{equation}
where $K_i\geq 0$, $\chi_i\in\R$ with $i=1,\ldots,k-1$. Then, we can choose gains $K_i>0$ and constants $\chi_1^*<-1$, $\chi_j=0$, $j=2,\ldots,k-1$, such that for $\ve>0$ sufficiently small, the origin is rendered locally asymptotically stable, but its region of attraction is larger compared to the choice $K=0$.
\end{theorem}

\begin{pf}
We shall prove the result for the case $k$ even. The proof for $k$ odd follows from the fact that if $k$ is odd, trajectories escape from a small neighborhood of the critical manifold when passing through singularities of even degeneracy. By following similar arguments as in the proof of Theorem \ref{thm:main2}, we can show that the origin of the blown up vector field in the chart $K_{\bar \ve}$ is still locally asymptotically stable for $K\geq 0$. This is due to the fact that $w$ is of order $O(\br)$ as $\br\to 0$. Next, note in \eqref{eq:th2} that for trajectories that do not converge to the origin, the term $-\frac{1}{\ve}z^k$ dominates the vector field $\dot z$. This means that such trajectories diverge from a small neighborhood of the origin with $z\to -\infty$. So, we look at the chart $K_{-\bar{z}}$. Thus, let us define the chart-coordinates 
\renewcommand{\bx}{\chi}
\renewcommand{\be}{\mu}
\renewcommand{\br}{\rho}
\begin{equation}\label{eq:buz}
(x_1,\ldots,x_{k-1},z,\ve)=(\br^k\bx_1,\ldots,\br^2\bx_{k-1},-\br,\br^{2k-1}\be),
\end{equation}
for $z<0$. Then, the corresponding blown up vector field reads as
\begin{equation}\label{eq:kz}
\begin{split}
\br' &= \br F(\bx)\\
\bx_i' &= H_i+G_i-K\br^{i-1}\be\bar w_i\\
\be' &= -(2k-1)\be F(\bx)
\end{split}
\end{equation}%
where $\bar w_i=\bar w_i(\br,\bx,\be)$ is induced by the blow up map \eqref{eq:buz}, that is $\bar w_i=w_i(\br^k\bx_1,\ldots,\br^2\bx_{k-1},-\br,\br^{2k-1}\be)$, and
\begin{equation}
\begin{split}
H_i&=\br^{i-1}\be(\bar f_i-C-b\be^{-\frac{1}{2k-1}}\delta_{1i}-\br^{1-i}\be^{-\frac{k}{2k-1}}a_i\bx_i)\\
G_i&=-(k-i+1)F(\bx)\bx_i\\
F(\bx)&=1-\sum_{i=1}^{2k-1}(-1)^i\bx_i,
\end{split}
\end{equation}%
with $\delta_{1i}=1$ if $i=1$ and $\delta_{1i}=0$ otherwise. Now we have the following crucial observation.
\begin{remark}\label{rem:rz}
Suppose $\phi(t)=(\br(t),\bx_1(t),\ldots,\bx_{k-1}(t),\be(t))$ is a trajectory of \eqref{eq:kz} satisfying $\br(t)\to 0$ and $\bx_i(t)\to\bx_i^*$ with $|\bx_i^*|<\infty$ as $\tau\to\infty$ (here $\tau$ denotes the rescaled time of \eqref{eq:kz}). Then, due to the blow up \eqref{eq:buz}, the trajectory $\phi$ is equivalent to a trajectory of the SFCS \eqref{eq:th2} that converges to the origin as $t\to\infty$.
\end{remark}
It is straightforward to show that for $K=0$, there is indeed an non-empty set of initial conditions such that the corresponding trajectory is as described in Remark \ref{rem:rz} (after all, for $K=0$ we have an equivalent system to the one in the chart $K_{\bar\ve}$). So the task of $\bar w$ is to enlarge such a set. To design $\bar w_i$ we shall use "high-gain" arguments. Although more complicated controllers can be designed, we want to keep the arguments as simple as possible to showcase the technique rather than the design itself. Thus, note that the arguments of Remark \ref{rem:rz} are satisfied if, for example, $\bx_1^*<-1$ and $\bx_j^*=0$ for all $j=2,\ldots,k-1$. With this idea in mind, we propose $\bar w=-K\br^q(\bx-\bx^*)$, where $K>0$ is diagonal, $-\infty<\bx_1^*<-1$ and $\bx_j=0$ for all $j=2,\ldots,k-1$, and the integer $q>1$ shall be set later. Note then that for $\br>0$, $\be>0$ and some large $K>0$, the point $\bx^*$ is attracting and $\br\to 0$. To obtain the expression of $w$ in \eqref{eq:cw} we just blow down using the expressions $\br=-z$ and $\bx_i=x_i(-z)^{-k+i-1}$. Therefore
\begin{equation}
w_i=-K_i(-z)^q\left(x_i(-z)^{-k+i-1}-\bx_i^*\right).
\end{equation}%
In principle, any integer $q$ such that $q-k+i-1>0$ provides a smooth expression for the controller $w$. For simplicity we choose $q=k-i+2$ so that the final expression of the controller is
\begin{equation}
w_i=K_i\left(x_iz+(-z)^{k-i+2}\bx_i^*\right),
\end{equation}
which is as stated in \eqref{eq:cw}. \hfill\ensuremath\qed
\end{pf}

{\change{
\begin{remark}
To design the ``compensation'' $w$ in Theorem \ref{thm:reg} we chose a high-gain controller. This was done to keep the arguments as simple as possible. Naturally, more elaborate and precise controllers may be designed, but the general idea stays the same: ``the design of controllers in the blow up space is rather simple''. We remark, however, that the high-gain nature of $w$ is to be expected. In some sense, the role of $w$ is to capture those trajectories that are quickly diverging and force them to return to the origin.
\end{remark}
}}

%% file: subfiles/ex.tex
\subsection{Example 1}

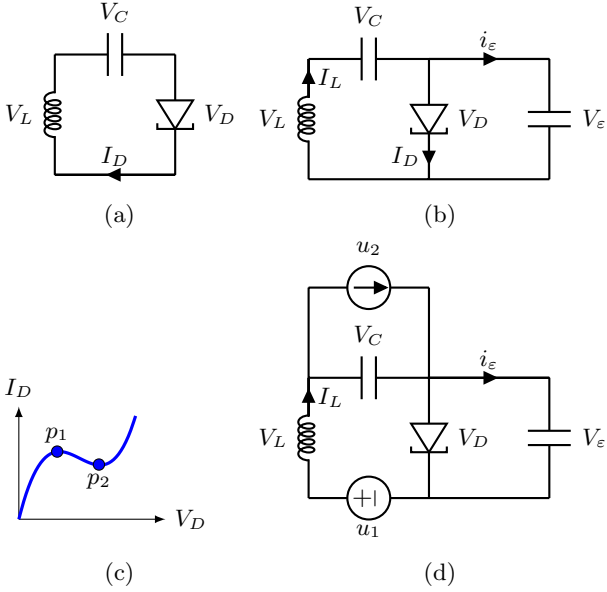
\begin{figure}[htbp]\centering
\begin{subfigure}[t]{0.39\columnwidth}
  \begin{tikzpicture}[scale=0.8]
    \draw[color=black,thick]
    (-1,-.5)
    to [L,l=$V_L$,bipoles/length=0.95cm,style={thin}] (-1,1.5)
    to [C,l=$V_C$,bipoles/length=0.95cm,style={thin}] (1,1.5)
    to [tDo,l=$V_D$,bipoles/length=0.95cm,style={thin}] (1,-.5)
    to [short,i_=$I_D$,style={ultra thin}] (-1,-.5)
    ;
  \end{tikzpicture}
  \caption{}
  \label{f:diode}  
\end{subfigure}\hfill
\begin{subfigure}[t]{0.6\columnwidth}
  \begin{tikzpicture}[scale=0.8]
    \draw[color=black,thick]
    (-1,-.5)
    to [L,l=$V_L$,bipoles/length=.95cm,style={thin}, i_>=$I_L$] (-1,1.5)
    to [C,l=$V_C$,bipoles/length=0.95cm,style={thin}] (1,1.5)
    to [tDo,l=$V_D$,bipoles/length=0.95cm,style={thin},, i_>=$I_D$] (1,-.5)
    to [short,style={ultra thin}] (-1,-.5)
    ;
    \draw[color=black,thick]
    (1,1.5) 
    to [short,i=$i_\ve$,style={ultra thin}] (3,1.5)
    to [C,l=$V_\ve$,bipoles/length=0.95cm,style={thin}] (3,-.5)
    to [short,style={ultra thin}] (1,-.5)
    ;
  \end{tikzpicture}
   \caption{}
  \label{f:diode-reg}  
  \end{subfigure}\\
\begin{subfigure}[t]{0.39\columnwidth}
  \begin{tikzpicture}[scale=0.15]
    \draw[-latex] (0,0)--(13,0) node[right] {$V_D$};
     \draw[-latex] (0,0)--(0,10) node[above] {$I_D$};
     \draw[scale=0.3,domain=0:5.75,smooth,variable=\x,blue,very thick] plot ({6*\x},{\x*\x*\x-9*\x*\x+24*\x});
     \draw[fill=blue] (3.4,6) circle (.5) node[above] {$p_1$};
     \draw[fill=blue] (7.1,4.85) circle (.5) node[below] {$p_2$};
  \end{tikzpicture}
  \caption{}
  \label{f:dgraph} 
\end{subfigure}
\begin{subfigure}[t]{0.6\columnwidth}
\begin{tikzpicture}[scale=0.8]
    \draw[color=black,thick]
    (-1,-.5)
    to [L,l=$V_L$,bipoles/length=.95cm,style={thin}, i_>=$I_L$] (-1,1.5)
    to [C,l=$V_C$,bipoles/length=.95cm,style={thin}] (1,1.5)
    to [tDo,l=$V_D$,bipoles/length=.95cm,style={thin}] (1,-.5)
    to [voltage source, v^<=$u_1$,bipoles/length=.95cm,style={thin}] (-1,-.5)
    (1,1.5) 
    to [short,i=$i_\ve$,style={ultra thin}] (3,1.5)
    to [C,l=$V_\ve$,bipoles/length=.95cm,style={thin}] (3,-.5)
    to [short,style={ultra thin}] (1,-.5)
    (-1,1.5)
    to [short,style={ultra thin}] (-1,3)
    to [current source,i=$u_2$,bipoles/length=.95cm,style={ thin}] (1,3)
    to [short,style={ultra thin}] (1,1.5)
    ;
  \end{tikzpicture}
  \caption{}
  \label{f:ex2}
\end{subfigure}
  \caption{(a) An electric circuit having a diode tunnel for negative resistance. (b) The characteristic curve of the tunnel diode. (c) Regularization of circuit (a). (d)The controlled circuit with voltage ($u_1$) and current ($u_2$) inputs.}
  \label{f:diode-c}
\end{figure}

We now exemplify the result of Theorem \ref{thm:main2} with an electric circuit having a tunnel diode \cite{Takens1,Reissig1} as shown in Figure \ref{f:diode}. 

\begin{figure}[htbp]\centering
\includegraphics[scale=0.95]{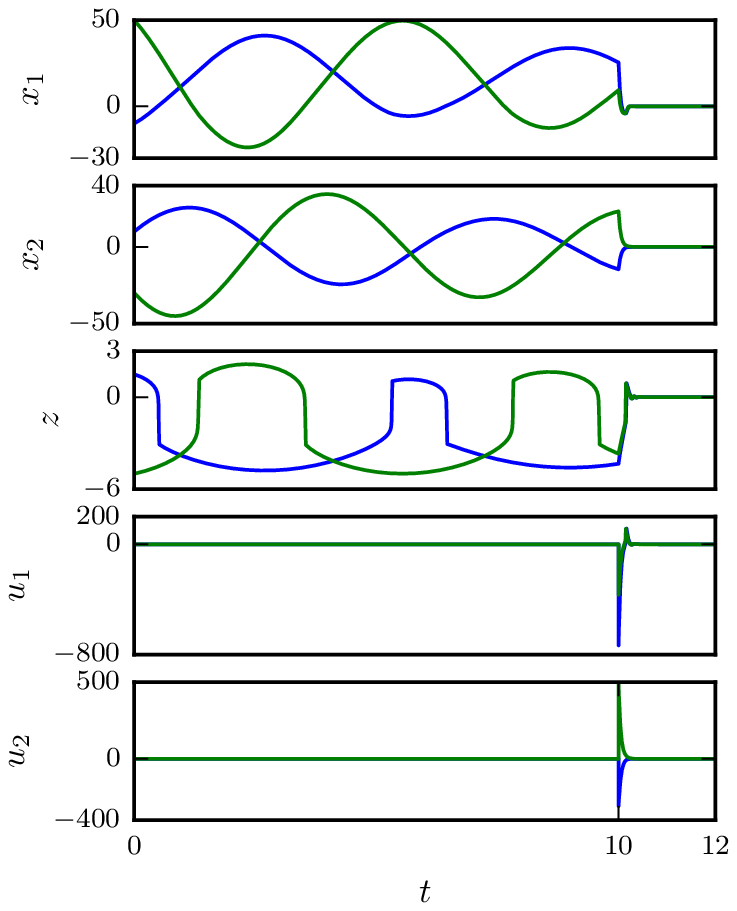}
\includegraphics[scale=0.95]{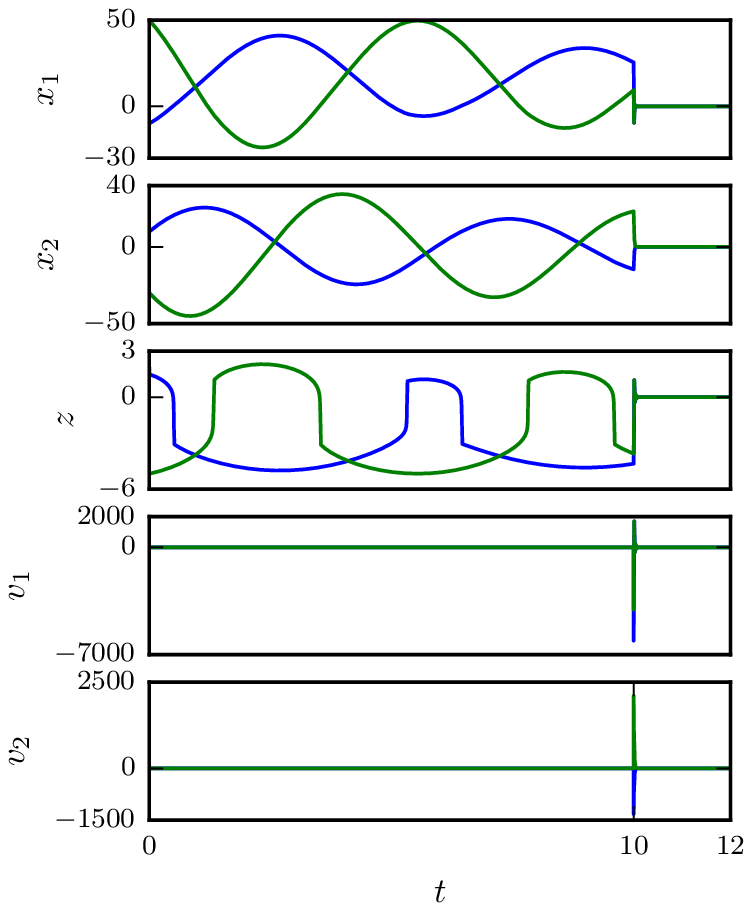}
\caption{Top: simulation of \eqref{eq:ex13} with the controllers \eqref{ex:c}. Bottom: simulation of \eqref{eq:ex13} with the high-gain controllers \eqref{ex:w}. 
We show trajectories with initial conditions $(x_1,x_2,z)=(-10,10,10)$ in blue and $(x_1,x_2,z)=(50,-30,-6)$ in green.
Note that, with a similar performance, $\max\left\{|u_i|\right\}<0.15\max\left\{|v_i|\right\}$, showing that the proposed controller in Theorem \ref{thm:main2} is considerably more convenient for implementation than a high-gain controller.   }
\label{fig:ex}
\end{figure}

In \cite{Reissig1} the diode's constitutive relation is given by $I_D=V_D^3-9V_D^2+24V_D$, see Figure \ref{f:dgraph}. A parasitic capacitance is added to regularize the circuit, as shown in Figure \ref{f:diode-reg}, see the justification in \cite{smale1972mathematical,ihrig1975regularization}. It is assumed that the parasitic capacitance $\ve$ is much smaller than any other parameter of the circuit. The equations describing the regularized circuit are
\begin{equation}\label{eq:ex11}
	\begin{split}
		\dot V_C &= \frac{1}{C}I_L\\
		\dot I_L &= -\frac{1}{L}(V_C+V_D)\\
		\ve\dot V_D &= -(V_D^3-9V_D^2+24V_D-I_L),
	\end{split}
\end{equation}
where $V_D$ denotes the voltage across the tunnel diode, $I_L$ is the current through the inductor, and $V_C$ the voltage across the capacitor. It is straightforward to show that \eqref{eq:ex11} has a unique equilibrium point at $p=(0,0,0)$, which is asymptotically stable. The critical manifold of \eqref{eq:ex11} is precisely given by the constitutive relation $$S=\left\{ I_L=I_D=V_D^3-9V_D^2+24V_D \right\}.$$ Note that $S$ has two fold points $p_1=(V_D,I_D)=(2,20)$ and $p_2=(V_D,I_D)=(4,16)$. The goal is to design a controller that stabilizes the operating point $(V_C,I_L,V_D)$ at one of the fold points, say $p_2$\footnote{Compare with \cite{GarciaCanseco2010127} where a similar {\change{diode system}} is studied. In there the authors stabilize a \emph{hyperbolic} operating point.} . The desired value of $V_C$ can be chosen arbitrarily but for simplicity we set it to  $V_C=0$. The controls are given by a voltage source ($u_1$) and acurrent source ($u_2$) as shown in Figure \ref{f:ex2}. 
Accordingly, the controlled system is described by
\begin{equation}\label{eq:ex12}
	\begin{split}
		\dot V_C &= \frac{1}{C}I_L +\frac{1}{C}u_2\\
		\dot I_L &= -\frac{1}{L}(V_C+V_D)+\frac{1}{L}u_1\\
		\ve\dot V_D &= -(V_D^3-9V_D^2+24V_D-I_L).
	\end{split}
\end{equation}%
For the analysis, let us perform the following change of coordinates $(x_1,x_2,z)=(V_C,-I_L+16,V_D-4)$. Thus, the operating point $p_2=(V_C,V_D,I_L)=(0,4,16)$ is translated to $(x_1,x_2,z)=(0,0,0)$. One then obtains
\begin{equation}\label{eq:ex13}
	\begin{split}
		\dot x_1&=\frac{1}{L}(x_2+z+4)+\frac{1}{L}u_1\\
		\dot x_2&=\frac{1}{C}(16-x_1)-\frac{1}{C}u_2\\
		\ve \dot z &= -(3z^2+x_1+z^3),
	\end{split}
\end{equation}
which is locally, near the origin, of the form studied in this paper. Using the results of Theorem \ref{thm:main2}, let us choose the controllers $u_1$ and $u_2$ as
\begin{equation}\label{ex:c}
\begin{split}
u_1 &= -4-\ve^{-\frac{2}{3}} a_1x_1+b\ve^{-\frac{1}{3}}z,\\
u_2&=16+\ve^{-\frac{2}{3}} a_2x_2.
\end{split}
\end{equation}
As a benchmark, we compare the performance of \eqref{ex:c} with high-gain controllers of the form
\begin{equation}\label{ex:w}
\begin{split}
\ve v_1 &= -A_1x_1+Bz,\\
\ve v_2&= -A_2x_2.
\end{split}
\end{equation}%
For the simulation shown in Figure \ref{fig:ex} we have chosen the parameters: $L=C=a_1=a_2=A_1=A_2=1$, $b=B=10$, $\ve=0.01$. We show trajectories for two initial conditions: $(x_1,x_2,z)=(-10,10,10)$ in blue and $(x_1,x_2,z)=(50,-30,-6)$ in green. We let the system run in open-loop for the first $10$ seconds. Then, at $t=10$ we ``turn on'' the controllers. We observe that both controllers ($u$ and $w$) provide a similar performance, however, the gains of $u_1$ and $u_2$ are approximately $10$ times smaller that those of $v_1$ and $v_2$.

\subsection{Example 2}

\begin{figure}[htbp]\flushleft
\begin{minipage}{.25\textwidth}
\includegraphics[scale=.9]{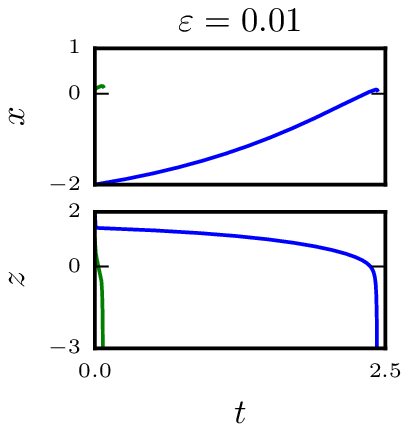}
\end{minipage}%
\begin{minipage}{.25\textwidth}
\includegraphics[scale=.9]{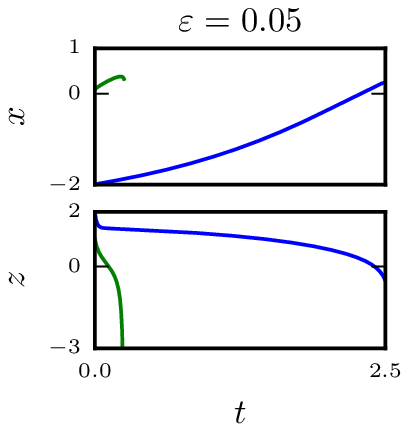}
\end{minipage}
\begin{minipage}{.25\textwidth}
\includegraphics[scale=.9]{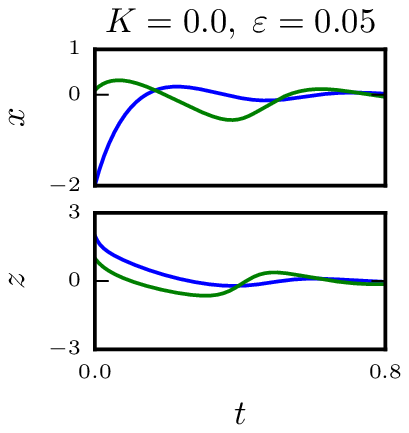}
\end{minipage}%
\begin{minipage}{.25\textwidth}
\includegraphics[scale=.9]{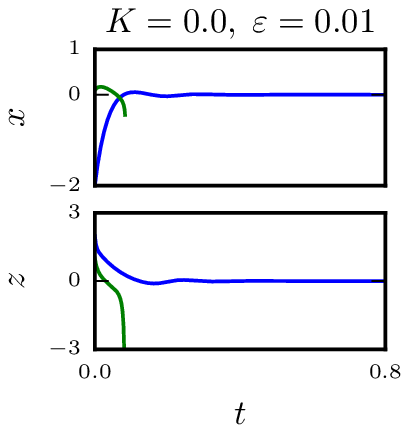}
\end{minipage}
\begin{minipage}{.25\textwidth}
\includegraphics[scale=.9]{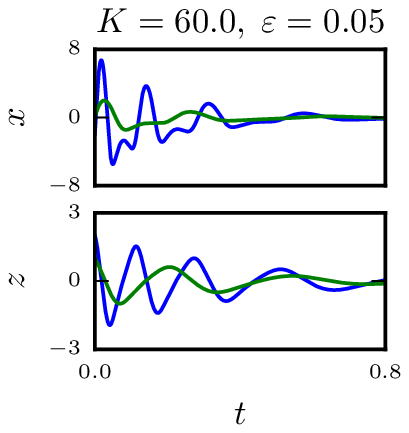}
\end{minipage}%
\begin{minipage}{.25\textwidth}
\includegraphics[scale=.9]{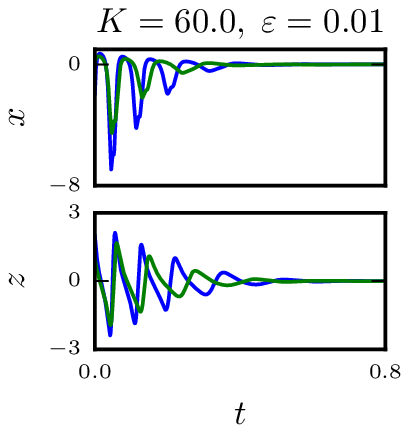}
\end{minipage}%
\caption{Simulation of the results of Theorem \ref{thm:reg} for $k=2$, and for initial conditions $(x,z)=(-2,2)$ (blue) and $(x,z)=(0.1,1)$ (green). The first row shows the open-loop dynamics. The second row shows the action of the controller of Theorem \ref{thm:reg} with $K=0$ (or equivalently of Theorem \ref{thm:main2}). The third row shows that the region of attraction of the origin has been enlarged due to the addition of $w$ as in \eqref{eq:cw}.} 
\label{fig:reg}
\end{figure}

To highlight the results of Theorem \ref{thm:reg}, let us consider the planar SFCS
\begin{equation}
\begin{split}
\dot x &= 1+x+z+u\\
\ve\dot z &= -(z^2+x)
\end{split}
\end{equation}
In Figure \ref{fig:reg}, we compare the performance of the controller \eqref{eq:cuu}-\eqref{eq:cw} with $K=0$ and $K>0$. For this simulation we  set the constants $A=1$, $b=3$ and $\chi^*=-2$. The first row of Figure \ref{fig:reg} shows the open-loop dynamics, and we observe that trajectories are quickly unbounded after crossing $z=0$. Next, in the second row of Figure \ref{fig:reg} we show the dynamics of the closed-loop system with controller proposed in Theorem \ref{thm:reg} with $K=0$. Note that for $\ve=0.05$, both trajectories converge to the origin, however, when we decrease $\ve$ to $\ve=0.01$, one trajectory diverges. Finally in the third row of Figure \ref{fig:reg} we show the effect of the compensation proposed in Theorem \ref{thm:reg}, and note that the origin is asymptotically stable for both values of $\ve$. In all these simulations we show trajectories with initial conditions $(x,z)=(-2,2)$ in blue and $(x,z)=(0.1,1)$ in green.

%% file: main.bbl
\begin{thebibliography}{10}

\bibitem{arnold1974critical}
V.~I. Arnold.
\newblock Critical points of smooth functions.
\newblock In {\em Proceedings of ICM-74}, volume~1, pages 19--40, 1974.

\bibitem{Arnold_singularities}
V.~I. Arnold, S.~M. Gusein-Zade, and A.~N. Varchenko.
\newblock {\em {Singularities of Differentiable Maps, Volume I}}, volume~17.
\newblock Birkh{\"{a}}user, 1985.

\bibitem{Brocker}
T.~Br{\"{o}}cker.
\newblock {\em {Differentiable Germs and Catastrophes}}, volume~17 of {\em
  Lecture Note Series}.
\newblock Cambridge University Press, 1975.

\bibitem{DumRou2}
F.~Dumortier and R.~Roussarie.
\newblock {\em {Canard Cycles and Center Manifolds}}, volume 121.
\newblock American Mathematical Society, 1996.

\bibitem{Fenichel1979}
N.~Fenichel.
\newblock {Geometric singular perturbation theory for ordinary differential
  equations}.
\newblock {\em Journal of Differential Equations}, 31(1):53--98, 1 1979.

\bibitem{GarciaCanseco2010127}
E.~Garc\'ia-Canseco, D.~Jeltsema, R.~Ortega, and J.~M.~A. Scherpen.
\newblock Power-based control of physical systems.
\newblock {\em Automatica}, 46(1):127 -- 132, 2010.

\bibitem{Hironaka1}
Heisuke Hironaka.
\newblock Resolution of singularities of an algebraic variety over a field of
  characteristic zero: I.
\newblock {\em Annals of Mathematics}, 79(1):109--203, 1964.

\bibitem{Hironaka2}
Heisuke Hironaka.
\newblock Resolution of singularities of an algebraic variety over a field of
  characteristic zero: Ii.
\newblock {\em Annals of Mathematics}, 79(2):205--326, 1964.

\bibitem{ihrig1975regularization}
E.~Ihrig.
\newblock The regularization of nonlinear electrical circuits.
\newblock {\em Proceedings of the American Mathematical Society}, pages
  179--183, 1975.

\bibitem{Jardon-Kojakhmetov2015}
H.~Jard{\'{o}}n-Kojakhmetov.
\newblock {Formal normal form of Ak slow-fast systems}.
\newblock {\em Comptes Rendus Mathematique}, 353(9):795--800, 9 2015.

\bibitem{Jardon-Kojakhmetov2016AnalysisSingularity}
H.~Jard{\'{o}}n-Kojakhmetov, H.~W. Broer, and R.~Roussarie.
\newblock {Analysis of a slow-fast system near a cusp singularity}.
\newblock {\em Journal of Differential Equations}, 260(4):3785--3843, 2016.

\bibitem{JardonMTNS2016}
H.~Jard{\'{o}}n-Kojakhmetov and J.~M.~A. Scherpen.
\newblock Stabilization of a planar slow-fast system at a non-hyperbolic point.
\newblock In {\em Proceedings of the 22nd International Symposium on
  Mathematical Theory of Networks and Systems}, July 2016.

\bibitem{jardon2017modelred}
H.~Jard\'on-Kojakhmetov and J.~M.~A. Scherpen.
\newblock {Model Order Reduction and Composite Control for a Class of Slow-Fast
  Systems Around a Non-Hyperbolic Point}.
\newblock {\em IEEE Control Systems Letters}, 1(1):68--73, July 2017.

\bibitem{JardonACC2017}
H.~Jard{\'{o}}n-Kojakhmetov, J.~M.~A. Scherpen, and D.~del Puerto-Flores.
\newblock Nonlinear adaptive stabilization of a class of planar slow-fast
  systems at a non-hyperbolic point.
\newblock In {\em Proceedings of the American Control Conference}, 2017.

\bibitem{jardon2017stabilization}
H.~{Jard\'on-Kojakhmetov}, J.~M.~A. {Scherpen}, and D.~{del Puerto-Flores}.
\newblock {Stabilization of slow-fast systems at fold points}.
\newblock {\em ArXiv e-prints}, April 2017.

\bibitem{KokotovicApps}
P.~V. Kokotovic.
\newblock {Applications of Singular Perturbation Techniques to Control
  Problems}.
\newblock {\em SIAM Review}, 26(4):501--550, 1984.

\bibitem{Kokotovic:1986:SPM:576779}
P.~V. Kokotovic, J.~O'Reilly, and H.~K. Khalil.
\newblock {\em {Singular Perturbation Methods in Control: Analysis and
  Design}}.
\newblock Academic Press, Inc., Orlando, FL, USA, 1986.

\bibitem{kosiuk2016geometric}
Ilona Kosiuk and Peter Szmolyan.
\newblock Geometric analysis of the goldbeter minimal model for the embryonic
  cell cycle.
\newblock {\em Journal of mathematical biology}, 72(5):1337--1368, 2016.

\bibitem{Krupa1}
M.~Krupa and P.~Szmolyan.
\newblock {Extending geometric singular perturbation theory to non hyperbolic
  points: fold and canard points in two dimensions}.
\newblock {\em SIAM J. Math. Anal.}, 33:286--314, 2001.

\bibitem{Krupa20102841}
M.~Krupa and M.~Wechselberger.
\newblock {Local analysis near a folded saddle-node singularity}.
\newblock {\em Journal of Differential Equations}, 248(12):2841--2888, 2010.

\bibitem{Kuehn2015}
C.~Kuehn.
\newblock {\em Multiple Time Scale Dynamics}.
\newblock Springer International Publishing, 2015.

\bibitem{Marino1988}
R.~Marino and P.~V. Kokotovic.
\newblock {A geometric approach to nonlinear singularly perturbed control
  systems}.
\newblock {\em Automatica}, 24(1):31--41, 1 1988.

\bibitem{Reissig1}
G.~Reissig.
\newblock Differential-algebraic equations and impasse points.
\newblock {\em IEEE Transactions on Circuits and Systems I: Fundamental Theory
  and Applications}, 43(2):122--133, Feb 1996.

\bibitem{RIAZA2011}
R.~Riaza.
\newblock {Explicit ODE reduction of memristive systems}.
\newblock {\em International Journal of Bifurcation and Chaos},
  21(03):917--930, 3 2011.

\bibitem{Rotstein2013}
H.~G. Rotstein.
\newblock {\em Mixed-Mode Oscillations in Single Neurons}, pages 1--9.
\newblock Springer New York, New York, NY, 2013.

\bibitem{Saksena1984}
V.~R. Saksena, J.~O'Reilly, and P.~V. Kokotovic.
\newblock {Singular perturbations and time-scale methods in control theory:
  Survey 1976-1983}.
\newblock {\em Automatica}, 20(3):273--293, 1984.

\bibitem{Shilnikov2012}
A.~Shilnikov.
\newblock Complete dynamical analysis of a neuron model.
\newblock {\em Nonlinear Dynamics}, 68(3):305--328, 2012.

\bibitem{smale1972mathematical}
S.~Smale.
\newblock On the mathematical foundations of electrical circuit theory.
\newblock {\em Journal of Differential Geometry}, 7(1-2):193--210, 1972.

\bibitem{Takens1}
F.~Takens.
\newblock {Constrained Equations: a Study of Implicit Differential Equations
  and their Discontinuous Solutions}.
\newblock In {\em Structural Stability, the Theory of Catastrophes, and
  Applications in the Sciences}, LNM 525, pages 134--234. Springer-Verlag,
  1976.

\bibitem{vanderpol_heart}
B.~van~der Pol and J.~van~der Mark.
\newblock {The heartbeat considered as a relaxation oscillation, and an
  electrical model of the heart}.
\newblock {\em The London, Edinburgh, and Dublin Philosophical Magazine and
  Journal of Science}, Ser.7,6:763--775, 1928.

\bibitem{Valery}
V.~D. Yurkevich.
\newblock {A unified approach to two-time scale control systems design: a
  tutorial}.
\newblock In {\em 2nd IASTED Int. Multi-Conference Automation, Control and
  Applications}, pages 314--319, 2005.

\end{thebibliography}
